\documentclass[10pt, reqno]{amsart}

\usepackage[T1]{fontenc}
\usepackage{lmodern,amsmath,amssymb,mathtools}
\usepackage[a4paper,margin=26mm]{geometry}
\usepackage[hidelinks,hypertexnames=false]{hyperref}

\newcommand{\R}{\mathbb R}
\newcommand{\B}{\mathbb B}
\newcommand{\eps}{\varepsilon}
\newcommand{\dd}{\,d}
\newcommand{\norm}[2]{\lVert #1\rVert_{#2}}

\allowdisplaybreaks[2]
\numberwithin{equation}{section}
\newtheorem{theorem}{Theorem}[section]
\newtheorem{proposition}[theorem]{Proposition}
\newtheorem{lemma}[theorem]{Lemma}
\newtheorem{corollary}[theorem]{Corollary}
\theoremstyle{definition}

\theoremstyle{remark}
\newtheorem{remark}[theorem]{Remark}

\allowdisplaybreaks[2]

\begin{document}

\title[Retarded interaction and non-scattering]{The retarded interaction and non-scattering
	for 3D Maxwell--Dirac system}

	\author{Yonggeun Cho}
\address{Department of Mathematics, and Institute of Pure and Applied Mathematics,
	Jeonbuk National University, Jeonju, 54896, Republic of Korea}
\email{changocho@jbnu.ac.kr}
\author{Kiyeon Lee}
\address{Department of Mathematics, University of Michigan, Ann Arbor, MI, USA}
\email{kiyeonl@umich.edu}

\begin{abstract}
We prove an a priori obstruction to linear scattering for the massive
Maxwell--Dirac system in three space dimensions. A global solution with
nonzero charge cannot scatter linearly in $L^2$ if its full and homogeneous
potentials are $O(t^{-1})$, and the homogeneous potential is
$o(t^{-1})$ in compact timelike cones. No nonlinear asymptotic expansion
is assumed in this non-scattering theorem. We also identify the exact
leading interaction when the nonlinear current approaches a free Dirac
current in $L^1$. Its coefficient is a positive quadratic form in the
velocity density and is bounded below by the square of the charge
divided by $4\pi$. The proof uses free Dirac asymptotics and an estimate
for retarded source errors after pairing with a smooth free current.
\end{abstract}

\thanks{2020{\it Mathematics Subject Classification.} 35Q41, 35B40, 35L05}

\maketitle

\section{Introduction and main results}\label{sec:introduction}

We consider the massive Maxwell--Dirac system on $\R^{1+3}$, with mass
and coupling constants normalized to one. In Hamiltonian form our
conventions are
\begin{equation}\label{eq:MD}
 \begin{split}
 (-i\partial_t+H)\psi&=V\psi,
 \qquad H=\alpha\cdot D+\beta,\quad D=-i\nabla,
 \quad V=A_0I_4+\mathbf A\cdot\alpha,\\
 (\partial_t^2-\Delta)A_0&=-\rho,
 \qquad (\partial_t^2-\Delta)\mathbf A=\mathbf j,
 \end{split}
\end{equation}
where $\psi$ takes values in $\mathbb C^4$, the potential
$A=(A_0,\mathbf A)$ is real, and
\begin{equation}\label{eq:current-definition}
 \rho=|\psi|^2,\qquad
 \mathbf j=(\psi^\dagger\alpha^1\psi,
             \psi^\dagger\alpha^2\psi,
             \psi^\dagger\alpha^3\psi).
\end{equation}
The Hermitian Dirac matrices satisfy
\[
 \alpha^j\alpha^k+\alpha^k\alpha^j=2\delta_{jk}I_4,
 \qquad \alpha^j\beta+\beta\alpha^j=0,
 \qquad \beta^2=I_4.
\]
The Lorenz condition is $-\partial_tA_0+\nabla\cdot\mathbf A=0$.
The proof below uses the reduced equations \eqref{eq:MD}; the Lorenz
constraints are relevant when specifying physical initial data.

The long-range behavior of this system is well established.
Flato, Simon, and Taflin \cite{FST} constructed modified wave operators
and studied asymptotic completeness. More recent analyses include the
physical-space approach of Herr, Ifrim, and Spitz \cite{HIS} and the
space-time resonance approach of Cho and Lee \cite{CL}.
The present paper concerns a specific quantitative mechanism: the
retarded field generated by a nonzero Dirac current produces a positive
interaction of order $t^{-1}$. We isolate an argument for this mechanism
that does not require a previously constructed nonlinear phase or a
pointwise asymptotic formula for the Maxwell potential. We make no
claim that the underlying long-range phenomenon is new.

The key distinction is between an asymptotic statement for the spinor
and one for its current. Multiplication by a common real phase changes
the spinor but leaves both quantities in \eqref{eq:current-definition}
unchanged. Consequently a nonlinear solution can have asymptotically
free currents while failing to scatter linearly as a spinor. The first
main theorem identifies the interaction in this situation. The second
uses a hypothetical linear scattering state to obtain a contradiction,
and therefore does not assume current asymptotics in its statement.

For a scalar or vector source $F$, let
\begin{equation}\label{eq:retarded-definition}
 \mathcal RF(t,x)=\frac1{4\pi}\int_0^t(t-s)
       \int_{\mathbb S^2}F(s,x+(t-s)\omega)\dd\omega\dd s.
\end{equation}
Thus the Maxwell equations have the decomposition
\begin{equation}\label{eq:wave-decomposition}
 A_0=A_0^{\rm hom}-\mathcal R\rho,
 \qquad \mathbf A=\mathbf A^{\rm hom}+\mathcal R\mathbf j,
\end{equation}
where $A^{\rm hom}$ is the free wave with the same Cauchy data as $A$.
We use the pointwise norm $|A|_*:=|A_0|+|\mathbf A|$.
Our field hypotheses are
\begin{equation}\label{eq:field-hypotheses}
 \begin{gathered}
 \sup_{t\ge0}(1+t)
       \bigl(\norm{|A(t)|_*}{L^\infty}
             +\norm{|A^{\rm hom}(t)|_*}{L^\infty}\bigr)\le M<\infty,\\
 \lim_{t\to\infty}\sup_{|v|\le\kappa}
       t|A^{\rm hom}(t,tv)|_*=0
       \quad\text{for every }0<\kappa<1.
 \end{gathered}
\end{equation}
These assumptions concern the full nonlinear potential and its actual
homogeneous component. They do not prescribe a favorable free
scattering state. Elementary sufficient conditions for the homogeneous
component, including the outgoing Coulomb tail, are proved below.

Write $U(t)=e^{-itH}$ and, for $f\in L^2(\R^3;\mathbb C^4)$, denote
the current of $U(t)f$ by $(\rho_f,\mathbf j_f)$. Let
$\B=\{v\in\R^3:|v|<1\}$. The velocity density $q_f\ge0$ is defined
explicitly in \eqref{eq:velocity-profile}; it satisfies
\begin{equation}\label{eq:velocity-mass-intro}
 \int_\B q_f(v)\dd v=\norm{f}{L^2}^2.
\end{equation}
Set, for $v,w\in\B$ and $v\ne w$,
\begin{equation}\label{eq:kernel-intro}
 \mathcal D(v,w)=
 \sqrt{(1-v\cdot w)^2-(1-|v|^2)(1-|w|^2)},\qquad
 K(v,w)=\frac{1-v\cdot w}{\mathcal D(v,w)}.
\end{equation}
For a nonnegative density $q$, its interaction energy is
\begin{equation}\label{eq:energy-definition}
 \mathcal E(q)=\frac1{4\pi}
       \iint_{\B\times\B}K(v,w)q(v)q(w)\dd v\dd w.
\end{equation}
Initially this expression is allowed to be infinite. Its finiteness
under the hypotheses of the following theorem is part of the proof.

\begin{theorem}[Interaction limit from current asymptotics]
\label{thm:current-interaction}
Let $(\psi,A)$ be a global solution of \eqref{eq:MD} satisfying
\eqref{eq:wave-decomposition} and \eqref{eq:field-hypotheses}.
Assume that charge is conserved,
$\norm{\psi(t)}{L^2}^2=Q$, and that for some $f\in L^2$,
\begin{equation}\label{eq:current-scattering}
 \norm{\rho(t)-\rho_f(t)}{L^1}
 +\norm{\mathbf j(t)-\mathbf j_f(t)}{L^1}
 \longrightarrow0\qquad(t\to\infty).
\end{equation}
Then $\norm{f}{L^2}^2=Q$, $\mathcal E(q_f)<\infty$, and
\begin{equation}\label{eq:main-interaction}
 \lim_{t\to\infty}t\mathcal I(t)=\mathcal E(q_f)
 \ge\frac{Q^2}{4\pi},\qquad
 \mathcal I(t):=-\int_{\R^3}
       (A_0\rho+\mathbf A\cdot\mathbf j)\dd x.
\end{equation}
If $Q>0$, the displayed lower bound for $\mathcal E(q_f)$ is strict,
and in particular
\begin{equation}\label{eq:eventual-lower}
 \mathcal I(t)\ge\frac{Q^2}{8\pi t}
 \quad\text{for all sufficiently large }t.
\end{equation}
\end{theorem}

\begin{theorem}[Obstruction to linear scattering]\label{thm:non-scattering}
Let $(\psi,A)$ be a global solution of \eqref{eq:MD} satisfying
\eqref{eq:wave-decomposition}, \eqref{eq:field-hypotheses}, and charge
conservation. Assume the Dirac equation holds in the $L^2$ Duhamel sense.
If there is an $f\in L^2$ such that
\[
 \norm{\psi(t)-U(t)f}{L^2}\longrightarrow0,
\]
then $Q=0$. In particular, every solution in this class with nonzero
charge fails to scatter linearly in $L^2$.
\end{theorem}

Sufficiently regular classical solutions with finite charge satisfy the
identities used in these statements. More generally, the proofs require
only the displayed Duhamel identities, local integrability of $V\psi$
in $L^2$, and charge conservation. The main theorems are a priori
statements about global solutions; they do not assert global existence
or derive \eqref{eq:field-hypotheses} for arbitrary initial data.
No smallness assumption is used once these hypotheses hold.

The proof has two approximation steps with different purposes.
First, the observation current is fixed to be the current of a free
spinor with smooth compact Fourier support. Such a current has
$L^\infty$ size $O(t^{-3})$. A retarded source error tending to zero in
$L^1$ then contributes $o(t^{-1})$ to the pairing. With this observation
state fixed, the source state can be approximated independently in
$L^2$. Only after this limit is justified do we approximate the
observation state. This order prevents constants in dispersive estimates
from being incorrectly treated as uniform over rough scattering states.

\section{Charge and elementary identities}\label{sec:identities}

For every unit vector $e\in\R^3$, $(e\cdot\alpha)^2=I_4$ and
$e\cdot\alpha$ is Hermitian. Hence
\begin{equation}\label{eq:causal-current}
 |\mathbf j|\le\rho.
\end{equation}
The Dirac equation and its adjoint give
\begin{equation}\label{eq:continuity}
 \partial_t\rho+\nabla\cdot\mathbf j=0.
\end{equation}
Integrating against expanding spatial cutoffs and using
\eqref{eq:causal-current} proves charge conservation for regular
finite-charge solutions. In particular,
$Q=\norm{\psi_0}{L^2}^2$.

The retarded wave operator is positive on scalar sources and is
self-adjoint in its spatial variables at every fixed pair of times.
If $S(\tau)=\sin(\tau|D|)/|D|$, then
\begin{equation}\label{eq:wave-elementary}
 \norm{S(\tau)b}{L^\infty}\le\tau\norm{b}{L^\infty},
 \qquad
 \int_{\R^3}S(\tau)b\dd x=\tau\int_{\R^3}b\dd x
 \quad(b\in L^1).
\end{equation}
Both formulas follow directly from its spherical-mean representation.
These simple estimates, applied after pairing, replace a pointwise
comparison of two retarded fields.

The sign of the interaction in \eqref{eq:main-interaction} follows
directly from \eqref{eq:wave-decomposition}. Its retarded part is
\[
 \int_{\R^3}\bigl(\rho\,\mathcal R\rho
               -\mathbf j\cdot\mathcal R\mathbf j\bigr)\dd x.
\]
It is nonnegative: $|\mathcal R\mathbf j|\le\mathcal R\rho$ and
\eqref{eq:causal-current} apply pointwise. A quantitative lower bound
requires more information, which is provided by the free-current
calculation in the next section.

There is also a scalar lower bound that requires no scattering
assumption and no decay estimate for the potential.

\begin{proposition}[Null-cone charge bound]\label{prop:cone-flux}
For a regular finite-charge solution and every $t>0$, $x\in\R^3$,
\begin{equation}\label{eq:cone-flux-bound}
 \mathcal R\rho(t,x)\ge\frac1{8\pi t}
       \int_{|y-x|<t}\rho(0,y)\dd y.
\end{equation}
In particular, for every $\kappa<1$,
\[
 \liminf_{t\to\infty}\inf_{|x|\le\kappa t}
       t\mathcal R\rho(t,x)\ge\frac Q{8\pi}.
\]
\end{proposition}

\begin{proof}
Integrate the continuity equation on the shrinking balls
$\{|y-x|<t-s\}$, $0\le s<t$. With
$n=(y-x)/|y-x|$, their charge satisfies
\[
 -\frac{d}{ds}\int_{|y-x|<t-s}\rho(s,y)\dd y
 =\int_{|y-x|=t-s}(\rho+n\cdot\mathbf j)(s,y)\dd S_y.
\]
The charge in the ball tends to zero as $s\uparrow t$ for a regular
solution. Integrating in $s$ and using $|\mathbf j|\le\rho$ gives
\[
 \int_{|y-x|<t}\rho(0,y)\dd y
 \le2\int_0^t\int_{|y-x|=t-s}\rho(s,y)\dd S_y\dd s
 \le8\pi t\,\mathcal R\rho(t,x).
\]
For $|x|\le\kappa t$, the initial ball contains
$\{|y|<(1-\kappa)t\}$. Charge integrability proves the final assertion.
\end{proof}

This controls $-A_0^{\rm ret}=\mathcal R\rho$. A homogeneous field can
alter the total scalar potential, and the spatial current can cancel
part of its contribution to the full interaction. The following
sections address that full contraction.

\section{Free Dirac currents and retarded pairings}\label{sec:free}

We use the unitary Fourier transform
\[
 \widehat f(\xi)=(2\pi)^{-3/2}\int_{\R^3}e^{-ix\cdot\xi}f(x)\dd x.
\]
For the free group $U(t)=e^{-itH}$ introduced above, let
\[
 \Pi_\theta(\xi)=\frac12\left(I+
       \theta\frac{\alpha\cdot\xi+\beta}{\langle\xi\rangle}\right),
 \qquad \theta\in\{+1,-1\}.
\]
The anticommutation relations imply that these are orthogonal
projections, $\Pi_++\Pi_-=I$, and
\begin{equation}\label{eq:free-projector-current}
 \Pi_\theta(\xi)\alpha^l\Pi_\theta(\xi)
 =\frac{\theta\xi_l}{\langle\xi\rangle}\Pi_\theta(\xi).
\end{equation}
Indeed, if $M=\alpha\cdot\xi+\beta$, then
$M\alpha^l+\alpha^lM=2\xi_lI$; multiplication on both sides by
$\Pi_\theta$ proves the identity.

For $v\in\B$, set $h(v)=\sqrt{1-|v|^2}$ and
$\xi_\theta(v)=\theta v/h(v)$. Define
\begin{equation}\label{eq:velocity-profile}
 \begin{split}
 u_{f,\theta}(v)
 &=e^{-3\pi i\theta/4}h(v)^{-5/2}
       \Pi_\theta(\xi_\theta(v))\widehat f(\xi_\theta(v)),\\
 q_f(v)&=\sum_{\theta=\pm1}|u_{f,\theta}(v)|^2.
 \end{split}
\end{equation}
These definitions make sense almost everywhere for every
$f\in L^2(\R^3;\mathbb C^4)$. The absolute Jacobian of
$v\mapsto\xi_\theta(v)$ is $h(v)^{-5}$. Consequently,
\begin{equation}\label{eq:velocity-charge}
 \sum_\theta\norm{u_{f,\theta}}{L^2(\B)}^2
 =\norm{f}{L^2}^2,
 \qquad \int_\B q_f(v)\dd v=\norm{f}{L^2}^2.
\end{equation}
Furthermore, \eqref{eq:free-projector-current} gives
\begin{equation}\label{eq:velocity-current}
 u_{f,\theta}(v)^\dagger\alpha^l u_{f,\theta}(v)
 =v_l|u_{f,\theta}(v)|^2.
\end{equation}

For any spinor $z$, write
$J[z]=(\rho[z],\mathbf j[z])=(z^\dagger z,z^\dagger\boldsymbol\alpha z)$,
and put $J_f(t)=J[U(t)f]$. Norms of four-component currents below
use $|(a_0,\mathbf a)|=|a_0|+|\mathbf a|$; all equivalent choices
give the same estimates.

\begin{lemma}\label{lem:current-continuity}
For $f,g\in L^2$, one has
\begin{align}
 \norm{q_f-q_g}{L^1(\B)}
 &\le (\norm{f}{L^2}+\norm{g}{L^2})\norm{f-g}{L^2},
       \label{eq:q-continuity}\\
 \sup_{t\in\R}\norm{J_f(t)-J_g(t)}{L^1}
 &\lesssim(\norm{f}{L^2}+\norm{g}{L^2})\norm{f-g}{L^2}.
       \label{eq:free-current-continuity}
\end{align}
The same current estimate holds for any two $L^2$ spinors at a fixed
time, without requiring them to solve the free equation.
\end{lemma}

\begin{proof}
For any bounded matrix $M$,
\[
 |z^\dagger Mz-w^\dagger Mw|
 \le\|M\|_{\mathrm{op}}(|z|+|w|)|z-w|.
\]
Cauchy--Schwarz proves the current estimate, using the unitarity of
$U(t)$. Apply the same argument in the direct sum of the two velocity
profile spaces and use \eqref{eq:velocity-charge} to obtain
\eqref{eq:q-continuity}.
\end{proof}

\begin{lemma}\label{lem:free-asymptotics}
Suppose $\widehat f\in C_c^\infty(\R^3;\mathbb C^4)$. Then
$u_{f,\theta}\in C_c^\infty(\B;\mathbb C^4)$, and, for $t\ge1$,
\begin{equation}\label{eq:free-spinor-expansion}
 U(t)f(x)=t^{-3/2}\sum_{\theta=\pm1}
 e^{-i\theta t h(x/t)}u_{f,\theta}(x/t)+r_f(t,x),
\end{equation}
where the leading term is defined to be zero outside $\B$ and
\begin{equation}\label{eq:free-spinor-remainder}
 \norm{r_f(t)}{L^\infty}\lesssim_f t^{-5/2},\qquad
 \norm{r_f(t)}{L^2}\lesssim_f t^{-1}.
\end{equation}
In particular, $\norm{U(t)f}{L^\infty}\lesssim_f t^{-3/2}$.
\end{lemma}

\begin{proof}
The $\theta$ half-wave is
\[
 (2\pi)^{-3/2}\int e^{it(v\cdot\xi-\theta\langle\xi\rangle)}
        \Pi_\theta(\xi)\widehat f(\xi)\dd\xi,
 \qquad v=x/t.
\]
Its stationary point is $\xi_\theta(v)$. The phase there is
$-\theta h(v)$, and its Hessian has signature $-3\theta$ and
absolute determinant $h(v)^5$. Stationary phase therefore gives
exactly the coefficient in \eqref{eq:velocity-profile}.

For clarity, only the following elementary, uniform version of
stationary phase is needed. On a compact set of velocities strictly
inside $\B$, write $\xi=\xi_\theta(v)+\eta$. Taylor's formula gives
the phase minus its critical value as
$-\theta\eta^TB_v(\eta)\eta/2$, where
\[
 B_v(\eta)=2\int_0^1(1-s)
       \nabla^2\langle\xi_\theta(v)+s\eta\rangle\dd s.
\]
This matrix is smooth and uniformly positive definite near $\eta=0$.
The map $z=B_v(\eta)^{1/2}\eta$ is therefore a smooth local
diffeomorphism, uniformly in these velocities. After localization the
phase becomes its critical value minus $\theta|z|^2/2$, with a smooth
compactly supported amplitude $b_v(z)$ whose relevant derivatives
are uniformly bounded. Fourier inversion of $b_v$ and the Gaussian
integral give
\[
 \int e^{-i\theta t|z|^2/2}b_v(z)\dd z
 =(2\pi/t)^{3/2}e^{-3\pi i\theta/4}
       \bigl(b_v(0)+O(t^{-1})\bigr).
\]
Indeed, the error is bounded by a constant times
$t^{-1}\int |\eta|^2|\widehat b_v(\eta)|\dd\eta$, using
$|e^{i\theta|\eta|^2/(2t)}-1|\le|\eta|^2/(2t)$.
On the remaining frequency region the phase gradient is separated
from zero, and repeated integration by parts gives arbitrarily
large inverse powers of $t$. This proves the $O(t^{-5/2})$
remainder uniformly on compact velocity sets.

Choose $\kappa<\kappa'<1$ so that every stationary velocity arising
from the Fourier support lies in $|v|\le\kappa$. If $|x|\ge\kappa't$,
the full phase gradient is bounded below by a constant times
$t+|x|$. Frequency integration by parts consequently gives
\[
 |U(t)f(x)|\lesssim_{f,N}(t+|x|)^{-N}.
\]
The leading term vanishes in this region. Combining this bound
with the uniform interior remainder proves the global
$L^\infty$ estimate. The interior region has volume $O(t^3)$,
so integration of the squared remainder gives the $L^2$ estimate;
the exterior estimate is stronger. The leading term has bounded
$L^2$ norm by \eqref{eq:velocity-charge}, which also verifies
all asserted integrability properties.
\end{proof}

\subsection{A retarded pairing identity}

The retarded operator \eqref{eq:retarded-definition} also has the form
\begin{equation}\label{eq:retarded-operator-free}
 \mathcal RF(t)=\int_0^t\frac{\sin((t-s)|D|)}{|D|}F(s)\dd s
 .
\end{equation}
By \eqref{eq:wave-elementary}, its spatial kernel has total mass
$t-s$ at source time $s$; hence
\begin{equation}\label{eq:retarded-L1}
 \norm{\mathcal RF(t)}{L^1}
 \le\int_0^t(t-s)\norm{F(s)}{L^1}\dd s.
\end{equation}
This also defines $\mathcal R$ on locally bounded
$L^1$-valued sources. In particular, all pairings below are
well-defined even when the source state belongs only to $L^2$.

Recall $\mathcal D$ and $K$ from \eqref{eq:kernel-intro}; their
values on the diagonal may be chosen arbitrarily. The identities
\begin{equation}\label{eq:kernel-comparison-free}
 \mathcal D(v,w)^2=|v-w|^2-|v\times w|^2
 \ge(1-|v|^2)|v-w|^2,
 \qquad K(v,w)\ge1
\end{equation}
hold off the diagonal. For the inequality, write
$v\times w=v\times(w-v)$. The last assertion follows from
$0<\mathcal D(v,w)\le1-v\cdot w$.

\begin{lemma}\label{lem:retarded-pairing}
Let $\widehat g\in C_c^\infty(\R^3;\mathbb C^4)$ and $f\in L^2$.
Then
\begin{equation}\label{eq:free-retarded-pairing-limit}
 \begin{split}
 B_t(g,f)&:=t\int_{\R^3}\bigl[
       \rho_g(t)\mathcal R\rho_f(t)
       -\mathbf j_g(t)\cdot\mathcal R\mathbf j_f(t)\bigr]\dd x,\\
 \lim_{t\to\infty}B_t(g,f)
 &=\frac1{4\pi}\iint_{\B\times\B}
       K(v,w)q_g(v)q_f(w)\dd v\dd w.
 \end{split}
\end{equation}
The integral is finite and is at least
$\norm{g}{L^2}^2\norm{f}{L^2}^2/(4\pi)$.
\end{lemma}

\begin{proof}
First suppose that $\widehat f$ is smooth and compactly supported.
Let $M_f$ be the current of the leading spinor in
\eqref{eq:free-spinor-expansion} for $s\ge1$, and put $M_f=0$
for $0\le s<1$. Write $J_f=M_f+E_f$, and define $M_g,E_g$
in the same way. The preceding lemmas give
\begin{equation}\label{eq:free-current-errors}
 \norm{E_f(s)}{L^1}\lesssim_f(1+s)^{-1},\qquad
 \norm{E_g(t)}{L^\infty}\lesssim_g t^{-4},\qquad
 \norm{M_g(t)}{L^\infty}\lesssim_g t^{-3}
 \quad(t\ge1).
\end{equation}
For example, the first estimate follows by Cauchy--Schwarz from
the $L^2$ spinor remainder, and the second from the
$L^\infty$ spinor remainder and $t^{-3/2}$ decay.

We may replace both currents in $B_t(g,f)$ by their leading
terms. In fact, $|\mathbf j_f|\le\rho_f$ and conservation of
the free $L^2$ norm imply
\[
 \int\mathcal R\rho_f(t,x)\dd x
 =\frac{t^2}{2}\norm{f}{L^2}^2.
\]
Thus the observation error, including the prefactor $t$ in
$B_t$, is $O_{f,g}(t^{-1})$. The source error is bounded by
\[
 Ct^{-2}\int_0^t(t-s)(1+s)^{-1}\dd s
 =O_{f,g}\bigl(t^{-1}\log(2+t)\bigr).
\]
These estimates use only \eqref{eq:retarded-L1}; no pointwise
asymptotic for a retarded potential has been assumed.

By \eqref{eq:velocity-current}, the leading current has the form
\begin{equation}\label{eq:leading-current-phases}
 M_f(s,sw)=s^{-3}\sum_{k\in\{0,2,-2\}}
       e^{iksh(w)}a_{f,k}(w),\qquad
 a_{f,0}(w)=(q_f(w),wq_f(w)),
\end{equation}
for $s\ge1$. All coefficients are smooth and compactly supported
in $\B$. Use the bilinear contraction
$\mathcal M(a,b)=a_0b_0-\mathbf a\cdot\mathbf b$.
The spatial kernel of the wave propagator is
$\delta(t-s-|x-y|)/(4\pi|x-y|)$. Set
$x=tv$, $y=sw$, and $s=t\sigma$. For $v\ne w$, the equation
\[
 1-\sigma=|v-\sigma w|,\qquad 0<\sigma<1,
\]
has the unique root
\begin{equation}\label{eq:retarded-root-free}
 \sigma(v,w)
 =\frac{1-v\cdot w-\mathcal D(v,w)}{1-|w|^2}
 =\frac{1-|v|^2}{1-v\cdot w+\mathcal D(v,w)}.
\end{equation}
The derivative of $1-\sigma-|v-\sigma w|$ at this root has
absolute value $\mathcal D(v,w)/|v-\sigma w|$. Consequently,
the integration in $\sigma$ gives the factor $1/\mathcal D(v,w)$.
On the compact velocity supports, $\sigma(v,w)$ is bounded below
by a positive constant, so the restriction $s\ge1$ disappears
for all sufficiently large $t$. The resulting integrals are
\begin{equation}\label{eq:oscillatory-retarded-pairing}
 \frac1{4\pi}\iint
 \frac{\mathcal M(a_{g,l}(v),a_{f,k}(w))}{\mathcal D(v,w)}
 e^{it[l h(v)+k\sigma(v,w)h(w)]}\dd v\dd w,
 \qquad l,k\in\{0,2,-2\}.
\end{equation}
The diagonal singularity is locally integrable by
\eqref{eq:kernel-comparison-free}. The delta-function calculation
can therefore first be made away from $v=w$ and then extended
to these integrals by absolute integrability.

We next show that every term with $(l,k)\ne(0,0)$ tends to zero.
Suppose first that $k\ne0$. For fixed $v$, write $w=v+r\omega$,
$r>0$, and set
\[
 b=\sqrt{1-|v\times\omega|^2},\qquad
 N=1-|v|^2-rv\cdot\omega,\qquad
 S(v,w)=\sigma(v,w)h(w).
\]
Then $\mathcal D(v,w)=rb$, and
\begin{equation}\label{eq:retarded-phase-monotonicity}
 S(v,v+r\omega)=\frac{N-rb}{h(v+r\omega)},\qquad
 \partial_rS(v,v+r\omega)
 =-\frac{bS(v,v+r\omega)}{1-|v+r\omega|^2}<0.
\end{equation}
Both formulas follow directly from \eqref{eq:retarded-root-free};
positivity follows from $N>\mathcal D$ inside $\B$.
Thus the change of variable $z=S(v,v+r\omega)$ is strictly
monotone on each ray. The radial amplitude, including the
Jacobian $r^2$ and the factor $1/\mathcal D$, is integrable.
Its pushforward is an $L^1$ function of $z$. The
Riemann--Lebesgue lemma makes its Fourier transform at frequency
$kt$ tend to zero. The factor $e^{itlh(v)}$ has modulus one;
dominated convergence in $v,\omega$ now proves the claim.

If $k=0$ and $l\ne0$, integrate the amplitude in $w$ first.
This leaves an $L^1$ function of $v$. Polar coordinates in $v$
and the substitution $z=\sqrt{1-|v|^2}$ again give the Fourier
transform of an $L^1$ function. This term also tends to zero.
No pointwise convergence of the mixed currents is used here.
For $l=k=0$, the contraction in
\eqref{eq:oscillatory-retarded-pairing} is
$(1-v\cdot w)q_g(v)q_f(w)$, proving the formula for smooth states.

We finally remove the regularity assumption on the source state.
For any current-valued source $E$, the decay of the fixed smooth
observation gives
\begin{equation}\label{eq:retarded-source-stability}
 t\left|\int\mathcal M(J_g(t),\mathcal RE(t))\dd x\right|
 \lesssim_g t^{-2}\int_0^t(t-s)\norm{E(s)}{L^1}\dd s,
 \qquad t\ge1.
\end{equation}
In particular, \eqref{eq:free-current-continuity} implies
\[
 \sup_{t\ge1}|B_t(g,f)-B_t(g,f_n)|
 \lesssim_g(\norm{f}{L^2}+\norm{f_n}{L^2})\norm{f-f_n}{L^2}.
\]
Choose $\widehat f_n\in C_c^\infty$ with $f_n\to f$ in $L^2$.
To pass to the limit on the right side of
\eqref{eq:free-retarded-pairing-limit}, put
\[
 K_g(w)=\int_\B K(v,w)q_g(v)\dd v.
\]
If $\operatorname{supp}q_g\subset\{|v|\le\kappa\}$, then
\eqref{eq:kernel-comparison-free} gives
\[
 0\le K_g(w)
 \le C_\kappa\int\frac{q_g(v)}{|v-w|}\dd v
 \lesssim_\kappa\norm{q_g}{L^\infty}+\norm{q_g}{L^1},
 \qquad w\in\B.
\]
Thus \eqref{eq:q-continuity} proves convergence of the limiting
integrals. This establishes the assertion for all $f\in L^2$.
The stated lower bound follows from $K\ge1$ and
\eqref{eq:velocity-charge}.
\end{proof}

\begin{remark}\label{rem:retarded-source-smallness}
An immediate consequence of \eqref{eq:retarded-source-stability}
will be useful below. If $\sup_{s\ge0}\norm{E(s)}{L^1}<\infty$
and $\norm{E(s)}{L^1}\to0$, then its right side tends to zero.
Indeed, split the time integral at a fixed $T$: the part $s\le T$
is $O(T/t)$, whereas the remaining part is bounded by a constant
times $\sup_{s\ge T}\norm{E(s)}{L^1}$. This proves smallness
of the retarded source error after pairing with a fixed smooth
free current, without requiring a bound on its trace on a light cone.
\end{remark}

\section{Transfer from the currents to the nonlinear interaction}\label{sec:transfer}

Throughout this section, $(\psi,A)$ solves \eqref{eq:MD} and satisfies
\eqref{eq:field-hypotheses}. We write
\[
 J_u=(\rho_u,\mathbf j_u),\qquad
 \rho_u=u^\dagger u,\qquad
 \mathbf j_u=u^\dagger\alpha u,
 \qquad u_f(t)=U(t)f.
\]
For any two spinors $u,v\in L^2(\R^3;\mathbb C^4)$, the Clifford
relations imply
\begin{equation}\label{eq:current-continuity}
 \norm{\rho_u-\rho_v}{L^1}
 +\norm{\mathbf j_u-\mathbf j_v}{L^1}
 \le 2\bigl(\norm{u}{L^2}+\norm{v}{L^2}\bigr)
             \norm{u-v}{L^2}.
\end{equation}
Indeed, $\norm{\omega\cdot\alpha}{\mathbb C^4\to\mathbb C^4}=1$
for $|\omega|=1$, and one expands each difference of quadratic forms
and applies Cauchy--Schwarz. We shall also use
$|\mathbf j_u|\le\rho_u$.

Suppose that \eqref{eq:current-scattering} holds with a free state $f$,
and put
\[
 d_f(s)=\norm{\rho_\psi(s)-\rho_{u_f}(s)}{L^1}
       +\norm{\mathbf j_\psi(s)-\mathbf j_{u_f}(s)}{L^1}.
\]
Charge conservation and \eqref{eq:current-scattering} give
\begin{equation}\label{eq:free-charge}
 \norm{f}{L^2}^2=Q,\qquad
 0\le d_f(s)\le4Q,\qquad d_f(s)\longrightarrow0.
\end{equation}
The bound on $d_f$ uses only the two conserved charges and
$|\mathbf j|\le\rho$.

Let
\[
 S(r)=\frac{\sin(r|D|)}{|D|},\qquad
 (\mathcal R F)(t)=\int_0^t S(t-s)F(s)\dd s.
\]
The signs in \eqref{eq:MD} give the exact decomposition
\begin{equation}\label{eq:field-decomposition}
 (A_0,\mathbf A)
 =(A_0^{\mathrm{hom}},\mathbf A^{\mathrm{hom}})
       +(-\mathcal R\rho_\psi,\mathcal R\mathbf j_\psi).
\end{equation}
No asymptotic expansion of the nonlinear spinor or potential is used
in this decomposition.

\begin{lemma}\label{lem:fixed-test-transfer}
For every $g$ with $\widehat g\in C_c^\infty(\R^3;\mathbb C^4)$,
\begin{equation}\label{eq:fixed-test-limit}
 \lim_{t\to\infty}
 -t\int_{\R^3}\bigl(A_0\rho_{u_g}
              +\mathbf A\cdot\mathbf j_{u_g}\bigr)(t,x)\dd x
 =\int_\B W_f(v)q_g(v)\dd v,
\end{equation}
where
\begin{equation}\label{eq:interaction-weight}
 W_f(v)=\frac1{4\pi}\int_\B
       \frac{1-v\cdot w}{\mathcal D(v,w)}q_f(w)\dd w.
\end{equation}
\end{lemma}

\begin{proof}
We first replace the source current in \eqref{eq:field-decomposition}.
In three dimensions, Kirchhoff's formula gives
\[
 (S(r)h)(x)=\frac{r}{4\pi}\int_{\mathbb S^2}h(x+r\omega)\dd\omega,
 \qquad \norm{S(r)h}{L^\infty}\le r\norm{h}{L^\infty}.
\]
The free stationary-phase estimate for the fixed smooth state $g$ is
$\norm{u_g(t)}{L^\infty}\le C_g\langle t\rangle^{-3/2}$.
Using the self-adjointness of $S(r)$, component by component, we obtain
\begin{align}\label{eq:source-error}
 &\left|\int_{\R^3}\left[
 \rho_{u_g}\mathcal R(\rho_\psi-\rho_{u_f})
 -\mathbf j_{u_g}\cdot\mathcal R(\mathbf j_\psi-\mathbf j_{u_f})
 \right](t,x)\dd x\right|\notag\\
 &\hspace{25mm}\le
 C_g\langle t\rangle^{-3}
          \int_0^t(t-s)d_f(s)\dd s=o(t^{-1}).
\end{align}
For the last equality, divide the time integral by $t^2$ and split it
at a fixed $T$. The part $s\le T$ tends to zero, whereas the remaining
part is bounded by $\frac12\sup_{s\ge T}d_f(s)$. Then let $T\to\infty$.
All the pairings above are absolutely convergent by the same estimate;
thus the interchange of the wave operator and the spatial pairing is
also justified for the $L^1$ current differences.

Next, the compact Fourier support of $g$ places all its group
velocities in a compact subset of $\B$. Consequently, for some
$\kappa<1$,
\begin{equation}\label{eq:smooth-free-tail}
 \int_{|x|>\kappa t}|u_g(t,x)|^2\dd x\longrightarrow0.
\end{equation}
This follows directly from nonstationary phase outside a slightly
larger compact set of group velocities. On $|x|\le\kappa t$, the
homogeneous-field hypothesis gives $tA^{\mathrm{hom}}=o(1)$ uniformly.
On its complement, use the global bound in
\eqref{eq:field-hypotheses} and \eqref{eq:smooth-free-tail}. It follows
that
\begin{equation}\label{eq:homogeneous-test-error}
 t\int_{\R^3}\bigl(A_0^{\mathrm{hom}}\rho_{u_g}
       +\mathbf A^{\mathrm{hom}}\cdot\mathbf j_{u_g}\bigr)\dd x
 \longrightarrow0.
\end{equation}
Equations \eqref{eq:field-decomposition}--\eqref{eq:homogeneous-test-error}
reduce the limit to the free retarded pairing $B_t(g,f)$.
Lemma~\ref{lem:retarded-pairing} identifies its limit with the
right-hand side of \eqref{eq:fixed-test-limit}.
\end{proof}

\begin{lemma}\label{lem:bounded-interaction-weight}
The function in \eqref{eq:interaction-weight} satisfies
\begin{equation}\label{eq:weight-bounds}
 \frac{Q}{4\pi}\le W_f(v)\le M
 \quad\text{for almost every }v\in\B.
\end{equation}
If $Q>0$, the first inequality is strict wherever $W_f$ is finite.
\end{lemma}

\begin{proof}
The kernel is nonnegative and satisfies
\[
 \frac{1-v\cdot w}{\mathcal D(v,w)}>1
 \qquad(v,w\in\B,\ v\ne w),
\]
because $(1-|v|^2)(1-|w|^2)>0$. Tonelli's theorem and
$\int_\B q_f\dd v=Q$ give the lower bound, including strictness
when $Q>0$. The value assigned to the kernel on the diagonal is
irrelevant for these integrals.

For the upper bound, \eqref{eq:field-hypotheses} and
$|\mathbf j_{u_g}|\le\rho_{u_g}$ give
\[
 \left|t\int_{\R^3}
       (A_0\rho_{u_g}+\mathbf A\cdot\mathbf j_{u_g})\dd x\right|
 \le M\norm{g}{L^2}^2.
\]
Taking the limit in Lemma~\ref{lem:fixed-test-transfer}, we obtain
\begin{equation}\label{eq:weight-test-bound}
 \int_\B W_f q_g\dd v\le M\int_\B q_g\dd v
 \qquad(\widehat g\in C_c^\infty).
\end{equation}
We explain why this tests all local scalar densities. Fix a sufficiently
small open set $V\Subset\B$ and one energy branch. The corresponding
rank-two spectral projector, evaluated at its momentum as a function
of velocity, has a smooth unit section $e(v)$ on $V$. Given
$h\in C_c^\infty(V)$, choose the Fourier data in that branch so that
its velocity profile is $h(v)e(v)$ and its other branch is zero.
The inverse velocity map and its positive Jacobian are smooth on $V$;
therefore these Fourier data belong to $C_c^\infty$ and satisfy
$q_g=|h|^2$. Thus \eqref{eq:weight-test-bound} implies
\[
 \int_V W_f(v)|h(v)|^2\dd v\le M\int_V|h(v)|^2\dd v.
\]
In particular $W_f$ is locally integrable, by taking a test function
equal to one on a smaller compact set. The last inequality then
implies $W_f\le M$ almost everywhere on $V$. A countable covering of
$\B$ proves the assertion.
\end{proof}

\begin{proof}[Proof of Theorem~\ref{thm:current-interaction}]
Set
\[
 \mathcal I(t)=-\int_{\R^3}
          (A_0\rho_\psi+\mathbf A\cdot\mathbf j_\psi)(t,x)\dd x.
\]
The current-scattering assumption and the global field bound give
\begin{equation}\label{eq:outer-current-replacement}
 t\mathcal I(t)
 =-t\int_{\R^3}(A_0\rho_{u_f}
                    +\mathbf A\cdot\mathbf j_{u_f})\dd x+o(1).
\end{equation}
Choose $g_n$ with $\widehat g_n\in C_c^\infty$ and $g_n\to f$ in $L^2$.
By \eqref{eq:current-continuity}, replacing $u_f$ by $u_{g_n}$ in the
right-hand side of \eqref{eq:outer-current-replacement} makes an error
bounded, uniformly in $t$, by
\[
 2M\bigl(\norm{f}{L^2}+\norm{g_n}{L^2}\bigr)
                      \norm{f-g_n}{L^2}.
\]
The velocity transformation and orthogonality of the spectral
projectors likewise give
\begin{equation}\label{eq:velocity-density-continuity}
 \norm{q_{g_n}-q_f}{L^1(\B)}
 \le\bigl(\norm{g_n}{L^2}+\norm{f}{L^2}\bigr)
                       \norm{g_n-f}{L^2}\longrightarrow0.
\end{equation}
The bound $W_f\le M$ therefore permits passage from $q_{g_n}$ to $q_f$
in \eqref{eq:fixed-test-limit}. First let $t\to\infty$ with $n$ fixed,
and then let $n\to\infty$. We conclude that
\begin{align}\label{eq:exact-interaction-limit}
 \lim_{t\to\infty}t\mathcal I(t)
 &=\int_\B W_f(v)q_f(v)\dd v\notag\\
 &=\frac1{4\pi}\iint_{\B\times\B}
       \frac{1-v\cdot w}{\mathcal D(v,w)}
             q_f(v)q_f(w)\dd v\dd w.
\end{align}
This is a finite number, bounded above by $MQ$ and below by
$Q^2/(4\pi)$. If $Q>0$, it is strictly larger than $Q^2/(4\pi)$.
This proves the claimed interaction limit and lower bound.
\end{proof}

\begin{proof}[Proof of Theorem~\ref{thm:non-scattering}]
Suppose, to the contrary, that
$\norm{\psi(t)-u_f(t)}{L^2}\to0$ for some $f\in L^2$.
Then \eqref{eq:current-continuity} proves
\eqref{eq:current-scattering}, and $\norm{f}{L^2}^2=Q$.
Define
\[
 \mathcal H(t)=\operatorname{Im}\int_{\R^3}u_f(t,x)^\dagger\psi(t,x)\dd x,
 \qquad V=A_0+\mathbf A\cdot\alpha.
\]
We have $|\mathcal H(t)|\le Q$. In the convention \eqref{eq:MD},
$\partial_t\psi=-iH\psi+iV\psi$, whereas
$\partial_tu_f=-iHu_f$. Self-adjointness of $H$ consequently gives
\[
 \mathcal H'(t)=\operatorname{Re}\int_{\R^3}u_f^\dagger V\psi\dd x
              =-\mathcal I(t)+o(t^{-1}).
\]
The error is bounded by
$\norm{u_f-\psi}{L^2}\norm{V}{L^\infty}Q^{1/2}=o(t^{-1})$.
For a general $L^2$ state $f$, the same identity follows from the
interaction representation and holds almost everywhere; equivalently,
one may first approximate $f$ by smooth states.
If $Q>0$, Theorem~\ref{thm:current-interaction} implies, for all
sufficiently large $t$,
\[
 \mathcal H'(t)\le-\frac{Q^2}{16\pi t}.
\]
Integration contradicts the boundedness of $\mathcal H$.
\end{proof}

\begin{corollary}\label{cor:common-phase-current}
Assume the hypotheses of Theorem~\ref{thm:current-interaction} except
\eqref{eq:current-scattering}. If, for some $f\in L^2$ and a real
measurable function $\Theta(t,x)$,
\[
 \norm{\psi(t)-e^{i\Theta(t,\cdot)}U(t)f}{L^2}\longrightarrow0,
\]
then the exact interaction limit \eqref{eq:exact-interaction-limit}
holds. The phase need not have any prescribed formula or regularity.
\end{corollary}

\begin{proof}
Multiplication by the same scalar phase in every spinor component
leaves both $\rho$ and $\mathbf j$ unchanged. Apply
\eqref{eq:current-continuity} to $\psi$ and $e^{i\Theta}U(t)f$, and then
Theorem~\ref{thm:current-interaction}.
\end{proof}

\begin{remark}
The argument proves an interaction and non-scattering criterion under
the field hypotheses \eqref{eq:field-hypotheses}. In particular, the
global $t^{-1}$ potential bound is an assumption here. No small-data
existence theorem, or derivation of that bound from initial data, is
asserted. The transfer uses convergence of the currents in $L^1$;
it requires neither nonlinear timelike profiles nor pointwise
convergence of those currents.
\end{remark}

\section{Conditions on the electromagnetic field}\label{sec:field-conditions}

We give explicit conditions on the free electromagnetic field that
enter Theorem~\ref{thm:non-scattering}. We also record a sufficient
condition on the nonlinear charge density for the required uniform
field bound. Neither statement uses a modified scattering theorem.
Throughout this section, $L=\partial_t^2-\Delta$.

\begin{lemma}[Decay of the free field]\label{lem:free-field-decay}
Let $W$ solve $LW=0$, with smooth initial data $(f,g)$. Suppose that,
for some $0<\sigma\le1$ and $M<\infty$,
\begin{equation}\label{eq:free-data-decay}
 |f(x)|\le M\langle x\rangle^{-1-\sigma},\qquad
 |\nabla f(x)|+|g(x)|\le M\langle x\rangle^{-2-\sigma}.
\end{equation}
Then
\begin{equation}\label{eq:free-global-decay}
 |W(t,x)|\lesssim_\sigma M(1+t+|x|)^{-1},\qquad t\ge0.
\end{equation}
Moreover, for every $0<\kappa<1$,
\begin{equation}\label{eq:free-timelike-decay}
 \sup_{|x|\le\kappa t}t|W(t,x)|
 \lesssim_{\sigma,\kappa}Mt^{-\sigma},\qquad t\ge1.
\end{equation}
The same conclusions hold componentwise for vector-valued data.
\end{lemma}

\begin{proof}
Kirchhoff's formula gives
\begin{equation}\label{eq:Kirchhoff-free}
 W(t,x)=\frac1{4\pi}\int_{\mathbb S^2}
 \bigl[f(x+t\omega)+t\omega\cdot\nabla f(x+t\omega)
                  +tg(x+t\omega)\bigr]\dd\omega.
\end{equation}
Write $r=|x|$. When $r\ge2t$ and $r\ge1$, the arguments of the
integrand have modulus at least $r/2$. The right-hand side is then
bounded by $CMr^{-1-\sigma}$. The case $t\ge2r$ and $t\ge1$ is
treated in the same way, with the bound $CMt^{-1-\sigma}$.
In the remaining region where $r$ and $t$ are comparable and
$r+t\ge2$, use
\[
 \int_{\mathbb S^2}\langle x+t\omega\rangle^{-p}\dd\omega
 =\frac{2\pi}{rt}
   \int_{|r-t|}^{r+t}s\langle s\rangle^{-p}\dd s.
\]
For $p=2+\sigma$, the last integral is bounded by $1/\sigma$.
Thus the two terms with a factor $t$ in \eqref{eq:Kirchhoff-free}
are bounded by $C_\sigma M/r$. For the remaining term, the
spherical mean identity
\[
 \int_{\mathbb S^2}|x+t\omega|^{-1}\dd\omega
 =\frac{4\pi}{\max(r,t)}
\]
and $\langle y\rangle^{-1-\sigma}\le |y|^{-1}$ give the same
bound. The bounded region $t+r\le2$ follows directly from
\eqref{eq:Kirchhoff-free}. This proves \eqref{eq:free-global-decay}.
If $|x|\le\kappa t$, every argument in that formula has modulus
at least $(1-\kappa)t$. Substitution of
\eqref{eq:free-data-decay} gives
$|W(t,x)|\lesssim_{\sigma,\kappa}Mt^{-1-\sigma}$, which proves
\eqref{eq:free-timelike-decay}.
\end{proof}

\subsection{The charge-compatible scalar data}
For the conventions
\[
 LA_0=-\rho,\qquad L\mathbf A=\mathbf j,\qquad
 -\partial_tA_0+\nabla\cdot\mathbf A=0,
\]
the initial constraints are
\begin{equation}\label{eq:field-Gauss}
 \dot a_0=\nabla\cdot\mathbf a,\qquad
 \Delta a_0=\nabla\cdot\dot{\mathbf a}+\rho_0.
\end{equation}
Indeed, these are the vanishing of the Lorenz expression and its
time derivative at $t=0$. The continuity equation
$\partial_t\rho+\nabla\cdot\mathbf j=0$ implies that the Lorenz
expression solves the homogeneous wave equation, so the constraints
propagate. We always take the Newtonian solution of the second
equation in \eqref{eq:field-Gauss}, with no harmonic addition.

Let $Q=\int_{\R^3}\rho_0\dd x$. Fix a smooth function $\chi$ with
$\chi(s)=0$ for $s\le1/2$ and $\chi(s)=1$ for $s\ge1$, and put
\begin{equation}\label{eq:outgoing-charge}
 C_Q(t,x)=-\frac{Q}{4\pi|x|}\chi(|x|-t),\qquad
 \widetilde A_\mu=A_\mu-\delta_{0\mu}C_Q.
\end{equation}
For $t\ge0$, $C_Q$ is a smooth free wave and vanishes when
$|x|\le t+1/2$. It also satisfies
$|C_Q(t,x)|\lesssim Q(1+t+|x|)^{-1}$.
Set
\[
 c(x)=-\frac{\chi(|x|)}{4\pi|x|},\qquad
 \varphi=\Delta c=-\frac{\chi''(|x|)}{4\pi|x|}.
\]
Then $\varphi$ is smooth and supported in an annulus, and
$\int\varphi\dd x=-\int_0^\infty r\chi''(r)\dd r=1$.
The corrected initial data are
\begin{equation}\label{eq:field-corrected-data}
 \begin{split}
 \widetilde a_0
 &=\Delta^{-1}\bigl(\nabla\cdot\dot{\mathbf a}
                           +\rho_0-Q\varphi\bigr),\\
 \dot{\widetilde a}_0
 &=\nabla\cdot\mathbf a-\frac{Q\chi'(|x|)}{4\pi|x|},
 \qquad
 (\widetilde{\mathbf a},\dot{\widetilde{\mathbf a}})
 =(\mathbf a,\dot{\mathbf a}).
 \end{split}
\end{equation}

If $\psi_0,\mathbf a,\dot{\mathbf a}$ are Schwartz, the source in
the first line of \eqref{eq:field-corrected-data} is Schwartz and
has integral zero. Its Newtonian potential is $O(\langle x\rangle^{-2})$
and its gradient is $O(\langle x\rangle^{-3})$.
To see this directly, subtract $|x|^{-1}$ from the Newtonian kernel
inside the convolution. On $|y|\le|x|/2$, the mean value theorem
bounds the resulting kernel and its $x$-gradient by
$C|y||x|^{-2}$ and $C|y||x|^{-3}$, respectively. The complementary
region has the same bounds by Schwartz decay and the local
integrability of $|x-y|^{-1}$ and $|x-y|^{-2}$.
All the other data in \eqref{eq:field-corrected-data} are Schwartz.
Consequently, Lemma~\ref{lem:free-field-decay}, with $\sigma=1$,
applies to the free wave generated by the corrected data. Adding
$\delta_{0\mu}C_Q$ shows that the physical homogeneous field satisfies
\eqref{eq:free-global-decay} and tends to zero after multiplication
by $t$, uniformly on every $|x|\le\kappa t$.

\subsection{A sufficient condition on the nonlinear density}
The next lemma gives one direct way to verify the global field bound.
Its assumption is a condition on the solution, not a consequence of
small initial data asserted in this paper.

\begin{lemma}[A retarded wave envelope]\label{lem:wave-envelope}
Let $U_R$ be the zero-data solution of
\[
 LU_R=\langle t\rangle^{-3}\mathbf1_{\{|x|\le t+R\}},
 \qquad R\ge1.
\]
Then
\begin{equation}\label{eq:radial-envelope-bound}
 \sup_{t\ge0}(1+t)\norm{U_R(t)}{L^\infty}
 \lesssim(1+R)^2.
\end{equation}
Consequently, if a Maxwell--Dirac solution satisfies
\begin{equation}\label{eq:density-envelope}
 \rho(t,x)\le M\langle t\rangle^{-3}
             \langle(|x|-t)_+\rangle^{-N},\qquad N>2,
\end{equation}
then its zero-data retarded field satisfies
\begin{equation}\label{eq:retarded-global-bound}
 \norm{A^{\mathrm{ret}}(t)}{L^\infty}
 \lesssim_N M(1+t)^{-1}.
\end{equation}
In particular, \eqref{eq:retarded-global-bound} applies if
$\psi_0$ has compact support and
$\norm{\psi(t)}{L^\infty}\lesssim\langle t\rangle^{-3/2}$.
\end{lemma}

\begin{proof}
Positivity of the retarded wave kernel gives $U_R\ge0$. Its radial
formula, for $r>0$, is
\[
 U_R(t,r)=\frac1{2r}\int_0^t\langle s\rangle^{-3}
  \int_{|r-t+s|}^{\min(r+t-s,s+R)}\lambda\dd\lambda\dd s,
\]
where an interval with upper endpoint below its lower endpoint
contributes zero. The value at $r=0$ is the continuous limit.
The elementary bound
$U_R(t,r)\le\int_0^t(t-s)\langle s\rangle^{-3}\dd s\lesssim t$
proves \eqref{eq:radial-envelope-bound} when $t\le4R$.
Suppose henceforth that $t>4R$.
If $r\le t/2$, the intersection with the source support forces
$s\ge(t-r-R)/2\ge t/8$, and therefore
\[
 U_R(t,r)\le\int_{t/8}^t(t-s)\langle s\rangle^{-3}\dd s
 \lesssim t^{-1}.
\]
If $r\ge t/2$, there is nothing to prove for $r>t+R$.
Otherwise set $T=t-r\ge-R$. The inner integral is bounded by
$(T+R)(s+R)$ and is zero unless $s\ge\max(0,(T-R)/2)$.
Thus
\[
 U_R(t,r)\lesssim\frac{T+R}{r}
  \int_{\max(0,(T-R)/2)}^\infty(s+R)\langle s\rangle^{-3}\dd s.
\]
For $T\le2R+2$, the product of $T+R$ and the last integral is
$O((1+R)^2)$. For $T>2R+2$, the lower endpoint is at least $T/4$;
the integral is $O(T^{-1}+RT^{-2})$, so that product is $O(1)$.
Since $r\ge t/2$, this proves \eqref{eq:radial-envelope-bound}.

For $z\ge0$,
\[
 \langle z\rangle^{-N}
 \le N\int_1^\infty R^{-N-1}\mathbf1_{\{z\le R\}}\dd R.
\]
By this inequality, positivity and
\eqref{eq:radial-envelope-bound} bound the retarded wave of
\eqref{eq:density-envelope} by
\[
 \frac{C_NM}{1+t}\int_1^\infty R^{-N-1}(1+R)^2\dd R
 \lesssim_N\frac{M}{1+t}.
\]
Since $|\mathbf j|\le\rho$, this estimate controls every component
of $A^{\mathrm{ret}}$, proving \eqref{eq:retarded-global-bound}.
Finally, the continuity equation gives
\[
 \frac{d}{dt}\int_{|x|>t+R}\rho(t,x)\dd x
 =\int_{|x|=t+R}(\mathbf j\cdot x/|x|-\rho)\dd S_x\le0.
\]
Thus the support of $\psi(t)$ lies in $|x|\le t+R$ if the initial
spinor is supported in $|x|\le R$. The stated dispersive bound supplies the
radial envelope used above.
\end{proof}

Combining Lemma~\ref{lem:wave-envelope} with
Lemma~\ref{lem:free-field-decay} and the outgoing charge wave gives
$\sup_{t\ge0}(1+t)\norm{A(t)}{L^\infty}<\infty$ whenever the
stated data and density conditions hold.

\section{An application using a known existence theorem}
\label{sec:external-application}

The criterion in Theorem~\ref{thm:non-scattering} and the preceding
field estimates are independent of a nonlinear asymptotic expansion.
The following application uses an external result only to obtain
global existence and the uniform decay bound for the full potential:
Theorem~1.1 and the charge-removal argument in Remark~1.5 of
Herr--Ifrim--Spitz~\cite{HIS}. Their modified scattering theorem is
not needed here.

\begin{corollary}\label{cor:localized-data-application}
Let $\psi_0\in\mathcal S(\R^3;\mathbb C^4)$ satisfy
\begin{equation}\label{eq:application-smallness}
 \norm{\langle x\rangle^3\psi_0}{H^{12}}\le\eps
\end{equation}
for sufficiently small $\eps>0$. Take
\[
 \mathbf a=\dot{\mathbf a}=0,\qquad
 a_0=\Delta^{-1}|\psi_0|^2,\qquad \dot a_0=0.
\]
Then the global solution supplied by the cited existence theorem
satisfies the hypotheses of Theorem~\ref{thm:non-scattering}.
In particular, if $\psi_0\ne0$, its spinor does not scatter in
$L^2$ to a free Dirac solution.
\end{corollary}

\begin{proof}
We verify the hypotheses of the external theorem after subtracting
\eqref{eq:outgoing-charge}. Write
\[
 \sigma_0=|\psi_0|^2-Q\varphi,
 \qquad h=\Delta^{-1}\sigma_0,
 \qquad b=-\frac{Q\chi'(|x|)}{4\pi|x|}.
\]
The corrected field data are $(\widetilde a_0,
\dot{\widetilde a}_0)=(h,b)$, with all spatial components zero.
We claim, for every fixed $0<\nu<1$,
\begin{equation}\label{eq:application-HIS-norms}
 \begin{split}
 \sum_{3j+k\le9}\Bigl(&
 \norm{|x|^j\partial^{\le j+k}\psi_0}{L^2}
 +\norm{|x|^j\partial^{j+k}h}{\dot H^{1/2}}
 +\norm{|x|^j\partial^{j+k}b}{\dot H^{-1/2}}\Bigr)\\
 &+\norm{h}{\dot H^{1/2-\nu}}
   +\norm{b}{\dot H^{-1/2-\nu}}
 \lesssim_\nu\eps+\eps^2.
 \end{split}
\end{equation}
Here derivatives are spatial and norms of derivative families are
summed. The spinor terms follow from
\eqref{eq:application-smallness} by the product rule. In particular,
all $x^\alpha\psi_0$, $|\alpha|\le3$, are bounded in $H^{12}$:
multiply $\langle x\rangle^3\psi_0$ by the smooth bounded
$H^{12}$ multiplier $x^\alpha/\langle x\rangle^3$.
Weighted product estimates give
\begin{equation}\label{eq:application-density-norms}
 \norm{\langle x\rangle^3\sigma_0}{L^1}
 +\sum_{|\alpha|\le3}\norm{x^\alpha\sigma_0}{H^8}
 \lesssim\eps^2,
 \qquad \int_{\R^3}\sigma_0\dd x=0.
\end{equation}
Consequently, on $|\xi|\le1$,
\[
 |\widehat{\sigma_0}(\xi)|\lesssim\eps^2|\xi|,
 \qquad
 |\partial_\xi^\gamma\widehat{\sigma_0}(\xi)|\lesssim\eps^2
 \quad(1\le|\gamma|\le3).
\]
For $|\alpha|=j$, $|\beta|=j+k$, expanding the derivatives in
\[
 |\xi|^{1/2}\partial_\xi^\alpha
       \bigl(\xi^\beta\widehat h(\xi)\bigr),\qquad
 \widehat h(\xi)=-|\xi|^{-2}\widehat{\sigma_0}(\xi),
\]
gives a bound $C\eps^2|\xi|^{k-1/2}$ at zero, which is square
integrable in three dimensions. At high frequencies, the coefficient
of $\partial_\xi^\gamma\widehat{\sigma_0}$ has degree
\[
 k+|\gamma|-\frac32\le k+j-\frac32
 \le\frac{15}{2}-2j\le\frac{15}{2}<8.
\]
Thus \eqref{eq:application-density-norms} proves the Cartesian
weighted norms of $h$. Its extra low-frequency norm is controlled
by the finite integral $\int_0^1r^{1-2\nu}\dd r$.
The function $b$ is smooth and supported in an annulus. All its
Cartesian weighted derivative norms are $O(Q)$, and its
$\dot H^{-1/2-\nu}$ norm is finite for the same reason, since
$\widehat b$ is bounded at zero.

For clarity, the Cartesian estimates also imply the literal radial
weights in \eqref{eq:application-HIS-norms}. A smooth degree-zero
angular multiplier is bounded on $\dot H^{1/2}(\R^3)$: the
Gagliardo seminorm bounds its commutator term by
$C\norm{|x|^{-1/2}f}{L^2}$, which is controlled by the fractional
Hardy inequality. Boundedness on $\dot H^{-1/2}$ follows by
duality. Now use
\[
 |x|^jf=\sum_{|\alpha|=j}\frac{j!}{\alpha!}
      \frac{x^\alpha}{|x|^j}\,x^\alpha f.
\]
This completes the proof of \eqref{eq:application-HIS-norms}.

These are exactly the smallness assumptions (1.8)--(1.9) of
\cite[Theorem~1.1]{HIS}. The outgoing charge correction is treated
by the argument of \cite[Remark~1.5]{HIS}; its sign is fixed by
\eqref{eq:field-Gauss} and \eqref{eq:outgoing-charge}.
To align the source signs when using that result, observe that
the printed reduced wave equation in \cite{HIS}, with its printed
definition of $\Box$, has the opposite source sign. The two
reduced systems are related by
\[
\psi^c=i\gamma^2\overline\psi,\qquad A^c=-A,\qquad
\gamma^2=\beta\alpha^2,
\]
in the standard Dirac representation, where the bar here denotes
entrywise complex conjugation. With $C=i\gamma^2$ the matrix identities
$C\overline{\alpha^j}C^{-1}=\alpha^j$ and
$C\overline\beta C^{-1}=-\beta$ show that this antiunitary transformation
anticommutes with the free Dirac Hamiltonian and preserves both
$\rho$ and $\mathbf j$. It therefore preserves the Dirac equation
with the potential changed to $-A$, while reversing the Maxwell
source sign. All data norms and absolute decay bounds are unchanged.
The charge-removal argument is used with the corresponding tail
sign; its estimates depend only on the magnitude of that tail.

The external theorem consequently gives global existence and
\[
 |A(t,x)|\lesssim\eps(1+t+|x|)^{-1}.
\]
The charge-compatible Schwartz data satisfy the homogeneous-field
conditions established after \eqref{eq:field-corrected-data}.
Theorem~\ref{thm:non-scattering} now applies. No nonlinear
timelike asymptotic expansion is used in this application.
\end{proof}

\end{document}